\documentclass[12pt]{amsart}
\usepackage{amssymb}
\usepackage{tabularx}
\usepackage{booktabs}

 \rm

\renewcommand{\(}{\left(}
\renewcommand{\)}{\right)}

\newcommand{\<}{\langle}
\renewcommand{\>}{\rangle}

\newcommand{\abs}[1]{\left\lvert#1\right\rvert}
\newcommand{\norm}[1]{\left\lVert#1\right\rVert}
\newcommand{\st}{\:|\:}

\renewcommand{\phi}{\varphi}

\newcommand{\C}{{\mathbb{C}}}

\newcommand{\N}{{\mathbb{N}}}
\theoremstyle{plain}
\newtheorem{thm}{Theorem}[section]
\newtheorem{lem}[thm]{Lemma}

\newtheorem{cor}[thm]{Corollary}

\theoremstyle{definition}

\theoremstyle{remark}

\title[Asymptotic Identities for Jacobi Polynomials]
{Asymptotic Identities for Jacobi Polynomials via Spectral Geometry of Rank-One Symmetric Spaces}

\author{Ankita~Sharma}
\email{ankitasharma.rs.mat19@itbhu.ac.in}
\thanks{The author was supported by a fellowship from CSIR
(File No.: 09/1217(0077)/2019-EMR-I)}
\address{Department of Mathematical Sciences, IIT(BHU), Varanasi, 221005,
INDIA}
\keywords{Compact rank-one symmetric space;
Laplace-Beltrami operator; Jacobi polynomial; Spherical function.}
\subjclass[2020]{33C50, 33C45, 33C55, 41A60}
\begin{document}

\begin{abstract}
Radial eigenfunctions of the Laplace-Beltrami operator on compact rank-one symmetric spaces
may be expressed in terms of Jacobi polynomials.
We use this fact to prove an identity for Jacobi polynomials
which is inspired by results of Minakshisundaram-Pleijel and Zelditch
on the Fourier coefficients of a smooth measure supported on
a compact submanifold of a compact Riemannian manifold.
%Also for the second approach we used the Christoffel-Darboux summation formula
%and Laplace's asymptotic formula for the Jacobi polynomials.

%asymptotic Formula
%for the Fourier coefficients of a smooth measure
%supported on a compact submanifold of a compact Riemannian manifold
%given by '' ''.
\end{abstract}

\maketitle

%\tableofcontents\newpage

\section{Introduction}
Jacobi polynomials are related to the {\em Laplace-Beltrami operator}
on compact rank-one symmetric spaces. In this paper, we use this relation to obtain some identities for the Jacobi polynomials.
\subsection{A brief review of Jacobi polynomials}\label{S:intro}
Given $\alpha>-1$, $\beta>-1$ and $\ell \in \N\cup\{0\}$, the Jacobi polynomial $\mathcal{P}_\ell^{(\alpha,\beta)}(x)$
may be defined by Rodrigues' formula
(see \cite[Equation (4.3.1)]{szego})
\begin{equation}\label{E:jacobi}
(1-x)^\alpha(1+x)^\beta \mathcal{P}_\ell^{(\alpha,\beta)}(x) =
\frac{(-1)^\ell}{2^\ell \ell!}\frac{d^\ell}{dx^\ell}\{(1-x)^{\ell+\alpha}(1+x)^{\ell+\beta}\}.
\end{equation}
The Jacobi polynomial $\mathcal{P}_\ell^{(\alpha,\beta)}(x)$ is a solution of the differential equation
(see \cite[Theorem 4.2.1]{szego})
\begin{equation}\label{E:jacobi de}
(1-x^2)y''+(\beta-\alpha-(\alpha+\beta+2)x)y'+\ell(\ell+\alpha+\beta+1)y=0.
\end{equation}
%Recall that the Jacobi operator is a differential operator. It is well known
%a positive selfadjoint second order linear differential operator
%in the weighted space $\mathcal{L}^2[(-1,1);(1-x)^\alpha(1+x)^\beta)dx]$.
%The spectrum of the Jacobi operator is given by the sequence of eigenvalues
%\begin{equation*}
%\tilde\lambda_\ell=\ell(\ell+\alpha+\beta+1), \; \ell>0.
%\end{equation*}
The Jacobi polynomials $\{\mathcal{P}_\ell^{(\alpha,\beta)}(x)\}_{\ell=0}^\infty$
are orthogonal on $[-1,1]$ with respect to the
weight function $(1-x)^\alpha(1+x)^\beta$ and satisfy the condition (see \cite[Equation (4.3.3)]{szego})
\begin{equation}\label{E:on}
\int_{-1}^1(1-x)^\alpha(1+x)^\beta \mathcal{P}_\ell^{(\alpha,\beta)}(x)^2~dx
=\frac{2^{\alpha+\beta+1}}{(2\ell+\alpha+\beta+1)}
\frac{\Gamma(\ell+\alpha+1)\Gamma(\ell+\beta+1)}{\Gamma(\ell+1)\Gamma(\ell+\alpha+\beta+1)}.
\end{equation}
For $\alpha>-1$ and $\beta>-1$, the Jacobi operator on
$\mathcal{L}^2[(-1,1);(1-x)^\alpha(1+x)^\beta)dx]$ is given by (see \cite[Equation (4.19)]{gangoli})
\begin{equation*}
J_{(\alpha,\beta)}=
(1-x^2)\frac{\partial^2}{\partial x^2}+(\beta-\alpha-(\alpha+\beta+2)x)\frac{\partial}{\partial x}.
\end{equation*}
%It is well known that the Jacobi operator is a non-negative operator on
%$\mathcal{L}^2[(-1,1);(1-x)^\alpha(1+x)^\beta)dx]$.
By Equation (\ref{E:jacobi de}), the spectrum of the Jacobi operator is given by the sequence of eigenvalues
\begin{equation*}
\tilde\lambda_\ell=-\ell(\ell+\alpha+\beta+1), \; \ell=0, 1, 2,\dots.
\end{equation*}

\subsection{Spherical functions on rank-one symmetric spaces in terms of Jacobi polynomials}\label{S:spherical}
Let $M$ be a compact rank-one symmetric space of real dimension $d$
with minimum sectional curvature $1$.
Then $M$ is homothetic (i.e, isometric up to a constant factor) to one of the following spaces
(see \cite{helgason(dgss)} and \cite[Theorem 8.4]{koranyi})
\begin{enumerate}
\item the $n$-sphere $\Sigma_n, \, n=1, 2, \dots$
\item the complex projective space $\C\mathbb{P}^n, \, n=2, 3, \dots$
\item the quaternionic projective space $\mathbb{H}\mathbb{P}^n, \, n=2, 3, \dots$
\item the Cayley projective plane $Ca\mathbb{P}^2$.
\end{enumerate}

%It is well known that the only compact rank-one symmetric spaces are
%the $n$-sphere $\Sigma_n$, the complex projective space $\C\mathbb{P}^n$,
%the quaternionic projective space $\mathbb{H}\mathbb{P}^n$ and
%the Cayley projective plane $Ca\mathbb{P}^2$,where $n \in \N$
%(see \cite{helgason(dgss)} and \cite[Theorem 8.4]{koranyi}).
Let $\rho$ be the distance function on $M$.
Let $\sigma$ be the Riemannian measure on $M$.
Fix a point $e \in M$. Recall that a radial function on $M$
is a function which depends only on $r=\rho(u,e)$ for $u \in M$.
Let $\Delta$ denote the Laplace-Beltrami operator on $M$.
Let $\{\lambda_\ell\}_{\ell=0}^\infty$ be the distinct eigenvalues of $-\Delta$, and let
$\mathcal{H}_\ell$ be the eigenspace corresponding to $\lambda_\ell$.
Let $m_\ell$ be the dimension of $\mathcal{H}_\ell$.
By the spectral theorem, the space $\mathcal{L}^2(M)$ is the
topological direct sum of the subspaces $\mathcal{H}_\ell$ (see e.g., \cite{roe}).
By \cite[Part 2 \S6, Proposition 2.10 and Corollary 3.3]{helgason(dgss)}
there exists a unique radial eigenfunction $\phi_\ell \in \mathcal{H}_\ell$
with $\norm{\phi_\ell}_{\mathcal{L}^2(M)}=1$; it
may be expressed in terms of the Jacobi polynomials,
as we explain below (for more details see \cite{gangoli}).
%Recall that the spherical functions on a rank-one symmetric space
%%($n$-sphere $\Sigma_n$, complex projective space $\C\mathbb{P}^n$,
%%quaternionic projective space $\mathbb{H}\mathbb{P}^n$ and the Cayley projective space $Ca\mathbb{P}^2$)
%are the radial eigenfunctions of the Laplace-Beltrami operator.
%The spherical functions on a compact rank-one symmetric space
%may be expressed in terms of the Jacobi polynomials (see \cite{gangoli}).

If $f$ is a radial $C^\infty$ function on $M$, then
%which depends only on $r=d(u,e)$ for $u \in M$
%the Laplace Beltrami operator of $M$ has the expression
(see \cite[Chapter \MakeUppercase{\romannumeral 10}, Lemma 7.12, \S 7]{helgason(dgss)}
and \cite[Equation (4.16)]{gangoli})
\begin{equation}\label{E:diff eq}
\Delta f(u) = \frac{\partial^2f}{\partial r^2} + \frac{A'(r)}{A(r)}\frac{\partial f}{\partial r},
\end{equation}
where $A(r)$ denotes the area of the sphere of radius $r$ centered at $e$ in $M$.
%Let $v \in T_oM$ be a unit tangent vector and
%$\gamma_v(r)$ the geodesic with $\gamma_v(0)=o$ and $\gamma'_v(0)=v$.
%Let $J(r)$ be the Riemannian density function (see \cite [\S3]{santhanam}) along $\gamma_v(r)$.
%(see \cite[Chapter \MakeUppercase{\romannumeral 1}, Lemma 4.10, \S4]{helgason})
Let $L$ be the diameter of $M$. Put $\omega = \frac{\pi}{2L}$.
Then for $0 \le r < L$
(cf. \cite{helgasonradon})
\begin{equation}\label{E:area}
A(r) = \frac{2\pi^\frac{p+q+1}{2}}{\Gamma(\frac{p+q+1}{2})\omega^{p+q}}\sin^{p+q}\omega r\cos^q\omega r,
\end{equation}
where $p$ and $q$ are non-negative integers depending on $M$(see Table \ref{table:parameters}).

\medskip

%\begin{center}
%\begin{table}[h!]
%\caption{Data for the symmetric spaces $M$(see \cite[page 171]{helgasonradon})}\label{table:parameters}
%\begin{tabular}{|c|c|c|c|c|c|}
%\hline
%%\toprule
%$M$                                       &$L$        &$d$     &$p$      &$q$      &$A'(r)/A(r)$\\
%%\toprule
%\hline
%$\Sigma_{n}$                        &$\pi$       &$n$     &$0$      &$d-1$   &$2\omega(n-1)\cot2\omega r$\\
%\hline
%$\C\mathbb{P}^n$                &$\pi/2$     &$2n$  &$d-2$   &$1$      &$\omega\((2n-1)\cot\omega r-\tan\omega r\)$\\
%\hline
%$\mathbb{H}\mathbb{P}^n$  &$\pi/2$     &$4n$  &$d-4$   &$3$      &$\omega\((4n-1)\cot\omega r-3\tan\omega r\)$\\
%\hline
%$Ca\mathbb{P}^2$                &$\pi/2$     &$16$  &$8$      &$7$     &$\omega\(15\cot\omega r-7\tan\omega r\)$\\
%\hline
%%\bottomrule
%\end{tabular}
%\end{table}
%\end{center}
\begin{center}
\begin{table}[h!]
\caption{Data for the symmetric spaces $M$(see \cite[page 171]{helgasonradon})}\label{table:parameters}
\begin{tabular}{|c|c|c|c|c|}
\hline
%\toprule
$M$                                       &$L$        &$d$     &$p$      &$q$   \\
%\toprule
\hline
$\Sigma_{n}$                        &$\pi$       &$n$     &$0$      &$d-1$ \\
\hline
$\C\mathbb{P}^n$                &$\pi/2$     &$2n$  &$d-2$   &$1$ \\
\hline
$\mathbb{H}\mathbb{P}^n$  &$\pi/2$     &$4n$  &$d-4$   &$3$ \\
\hline
$Ca\mathbb{P}^2$                &$\pi/2$     &$16$  &$8$      &$7$ \\
\hline
%\bottomrule
\end{tabular}
\end{table}
\end{center}

\medskip

%Therefore
%\begin{equation*}
%\frac{A'(r)}{A(r)}
%=\omega\((p+q)\cot\omega r - q\tan\omega r\).
%\end{equation*}

%Let $x=\cos r$ for $\Sigma_n$. Then Equation (\ref{E:diff eq}) reduces to
%\begin{equation}\label{E:de sphere}
%\Delta f = -(1-x^2)\frac{\partial^2f}{\partial x^2} +nx\frac{\partial f}{\partial x}.
%\end{equation}
%
%By comapring Equations (\ref{E:jacobi de}) and (\ref{E:de sphere}), we get for $\Sigma_n$
%\begin{equation*}
%\alpha=\frac{n-2}{2}=\beta.
%%\; \text{and} \; \tilde\lambda_\ell=\ell(\ell+n-1).
%\end{equation*}

If we take $x=\cos2\omega r$, then
%Then Equation (\ref{E:diff eq}) reduces to
\begin{equation}\label{E:de rank1}
\Delta f(u)= 4\omega^2J_{\(\frac{p+q-1}{2},\frac{q-1}{2}\)}f(x),
%\Delta f = -4\omega^2\[(1-x^2)\frac{\partial^2f}{\partial x^2} -\frac{1}{2}\(p+(p+2q+2)x\)\frac{\partial f}{\partial x}\].
\end{equation}
and so the function
$\mathcal{P}_\ell^{\(\frac{p+q-1}{2},\frac{q-1}{2}\)}(\cos2\omega\rho(u,e))$
is a radial eigenfunction of $\Delta$ with eigenvalue $4\omega^2\tilde\lambda_\ell$.
Therefore $\lambda_\ell=-4\omega^2\tilde\lambda_\ell$ and
$\phi_\ell(u)=c_\ell\mathcal{P}_\ell^{\(\frac{p+q-1}{2},\frac{q-1}{2}\)}(\cos2\omega\rho(u,e))$,
where $c_\ell$ is a normalizing constant.

\medskip

%\begin{center}
%\begin{table}[h!]
%\caption{Parameters related to $M$}
%\begin{tabular}{|c|c|c|c|}
%\hline
%%\toprule
%$M$                                               &  $\alpha$             & $\beta$\\
%%\toprule
%\hline
%$\Sigma_{n}$                                 &  $\frac{n-2}{2}$   & $\frac{n-2}{2}$\\
%\hline
%$\C\mathbb{P}^n$                          &  $n-1$                 & $0$\\
%\hline
%$\mathbb{H}\mathbb{P}^n$            &  $2n-1$               & $1$\\
%\hline
%$Ca\mathbb{P}^2$                          &  $7$                   & $3$\\
%\hline
%%\bottomrule
%\end{tabular}
%\end{table}
%\end{center}
%
%\medskip

\begin{lem}\label{L:projection}
The normalizing constant $c_\ell$ is given by
\begin{equation*}
c_\ell=
\sqrt{\frac{\omega^{p+q+1}\Gamma(\frac{p+q+1}{2})}{\pi^\frac{p+q+1}{2}}
\frac{(\frac{4\ell+p+2q}{2})\Gamma(\ell+1)\Gamma(\frac{2\ell+p+2q}{2})}
{\Gamma(\frac{2\ell+p+q+1}{2})\Gamma(\frac{2\ell+q+1}{2})}}.
\end{equation*}
\end{lem}
\begin{proof}
Let $\Pi:M \to [-1,1]$ be defined by
\begin{equation*}
\Pi(u)=\cos2\omega\rho(u,e).
\end{equation*}
Let $\Pi_1:M \to [0,L]$ and $\Pi_2:[0,L] \to [-1,1]$ be defined by
\begin{equation*}
\Pi_1(u)=\rho(u,e) \; \text{and} \; \Pi_2(r)=\cos{2\omega r}.
\end{equation*}
Then 
\begin{equation*}
\Pi={\Pi_2}\circ{\Pi_1}.
\end{equation*}
Therefore from Equation (\ref{E:area}), the pushforward of $\sigma$ by $\Pi$ is given by
\begin{equation*}
\begin{aligned} 
\Pi_*\sigma(x)=
&\;{\Pi_2}_*({\Pi_1}_*\sigma)\\
=&\;{\Pi_2}_*
\(\frac{2\pi^\frac{p+q+1}{2}}{\Gamma(\frac{p+q+1}{2})\omega^{p+q}}\sin^{p+q}\omega r\cos^q\omega r~dr\)\\
=&\;\frac{\pi^\frac{p+q+1}{2}(1-x)^\frac{p+q-1}{2}(1+x)^\frac{q-1}{2}}
{\Gamma(\frac{p+q+1}{2})
\omega^{p+q+1}2^{\frac{p+2q}{2}}}~dx.
\end{aligned}
\end{equation*}
It follows that
\begin{equation*}
\begin{aligned}
1=
&\;\norm{\phi_\ell}^2_{\mathcal{L}^2(M)}\\
=&\;\int_M c_\ell^2\abs{\mathcal{P}_\ell^{\(\frac{p+q-1}{2},\frac{q-1}{2}\)}(\cos2\omega\rho(u,e))}^2~d\sigma\\
=&\;\frac{c_\ell^2\pi^\frac{p+q+1}{2}}
{\Gamma(\frac{p+q+1}{2})\omega^{p+q+1}2^{\frac{p+2q}{2}}}
\int_{-1}^1(1-x)^\frac{p+q-1}{2}(1+x)^\frac{q-1}{2}
\abs{\mathcal{P}_\ell^{\(\frac{p+q-1}{2},\frac{q-1}{2}\)}(x)}^2~dx.
\end{aligned}
\end{equation*}
Using Equation (\ref{E:on}), we get
\begin{equation*}
c_\ell^2=
\frac{\omega^{p+q+1}\Gamma(\frac{p+q+1}{2})}{\pi^\frac{p+q+1}{2}}
\frac{(\frac{4\ell+p+2q}{2})\Gamma(\ell+1)\Gamma(\frac{2\ell+p+2q}{2})}
{\Gamma(\frac{2\ell+p+q+1}{2})\Gamma(\frac{2\ell+q+1}{2})}.
\end{equation*}
\end{proof}

\section{The identity}
Let $K$ be the group of all isometries of $M$ which fix $e$.
Let $dk$ be the Haar probability measure on $K$.
We need the following lemma to prove Theorem \ref{T:jacobimain}.
\begin{lem}\label{L:ip}
Let $\mu$ be a $K$-invariant measure on $M, \, i.e., \, k_*\mu=\mu$ for every $k \in K$.
Then for $\xi \in \mathcal{H}_\ell$,
\begin{equation*}
\<\xi, \phi_\ell\>=0 \implies \int_M\xi~d\mu=0.
\end{equation*}
\end{lem}
\begin{proof}
Let $\xi \in \mathcal{H}_\ell$ be such that $\<\xi, \phi_\ell\>=0$. Let
\begin{equation*}
\tilde\xi=\int_K\xi\circ k~dk.
\end{equation*}
Then $\tilde\xi$ is $K$-invariant. Therefore there exists a constant $a$ such that
$\tilde\xi=a\phi_\ell$.
Further,
\begin{equation*}
\begin{aligned}
\<\tilde{\xi}, \phi_\ell\>=
&\;\int_K\<\xi\circ k, \phi_\ell\>~dk\\
=&\;\int_K\<\xi, \phi_\ell \circ k^{-1}\>~dk\\
=&\;\int_K\<\xi, \phi_\ell\>~dk\\
=&\;0.
\end{aligned}
\end{equation*}
Therefore $\tilde\xi=0$.
Now for every $k \in K$,
\begin{equation*}
\begin{aligned}
\int_M\xi~d\mu=
&\;\int_K\int_M\xi~d\mu~dk\\
=&\;\int_K\int_M\xi~k_*d\mu~dk\\
=&\;\int_K\int_M\xi \circ k~d\mu~dk\\
=&\;\int_M\tilde\xi~d\mu\\
=&\;0.
\end{aligned}
\end{equation*}
\end{proof}

Let $Y_\ell^0=\phi_\ell$, and extend $\{Y_\ell^0\}$ to an orthonormal basis
$\{Y_\ell^i\}_{i=0}^{m_\ell-1}$ of $\mathcal{H}_\ell$. Let
\begin{equation*}
\{\Phi_j\}_{j=0}^\infty=\bigcup_{l=0}^\infty \{Y_\ell^i\}_{i=0}^{m_\ell-1}.
\end{equation*}

Then $\{\Phi_j\}_{j=0}^\infty$ is an orthonormal
basis of $\mathcal{L}^2(M)$ and
\begin{equation*}%\label{E:spectral-thm}
-\Delta \Phi_j = \gamma_j \Phi_j,
\end{equation*}
where $\gamma_j \in \{\lambda_i\}_{i=0}^\infty$.

%Thus, by Parseval's formula (see e.g., \cite[Theorem 4.18]{Rudin}),
%if $\psi\in C^\infty(M)$ then
%\begin{equation}\label{E:parseval}
%\sum_{\lambda_j < T}\abs{\widehat{\psi}(j)}^2 \sim
%\int_M \abs{\psi}^2 \,d\mu, \quad T\to\infty,
%\end{equation}
%where $\widehat{\psi}(j)=\<\psi,\Phi_j\>$ is the $j$-th Fourier
%coefficient of $\psi$.

Recall that if $\tau$ is a measure on $M$, the $j$-th Fourier coefficient
of $\tau$ (as a distribution on $M$) is (see \cite{sharma})
\begin{equation*}
\widehat{\tau}(j) = \<\tau, \Phi_j\> = \int_M \Phi_j \,d\tau.
\end{equation*}
Suppose $N$ is a compact submanifold of $M$ and let $\nu$ denote the Riemannian measure on $N$.
Let $\psi \in C^\infty(N)$. A measure of the form $\tau = \psi\nu$ is called a smooth measure
supported on $N$ (see e.g. \cite[Chapter 8, \S3]{stein}).
The following theorem is a result of Minakshisundaram-Pleijel and Zelditch (see \cite{sharma} and \cite{steven})
\begin{thm}\label{T:main}
Let $\tau=\psi\nu$ be a smooth measure supported on a compact
codimension $k$ submanifold $N$ of $M$.
Then
\begin{equation*}
\sum_{\gamma_j < T}\abs{\widehat{\tau}(j)}^2 \sim
\frac{T^{k/2}\int_N\abs{\psi}^2 d\nu}{(4\pi)^{k/2}\Gamma(\frac{k}{2}+1)},
\quad T\to \infty.
\end{equation*}
\end{thm}
%Let $\Pi:M \to [-1,1]$ be the projection map.
%Then push forward measure of $\mu$ with respect to
%$\Pi$ is given by
%\begin{equation*}
%\Pi_*\mu=\nu(N)dr.
%\end{equation*}
%Let $x=\cos r$ for $\Sigma_n$ and
%$x=\cos2r$ for $\C\mathbb{P}^n, \mathbb{H}\mathbb{P}^n$ and $Ca\mathbb{P}^2$.
%
%\begin{thm}\label{T:mainthm}
%For $\alpha>-1$ and $\beta>-1$,
%\begin{equation*}
%\lim_{m\to\infty} \frac{1}{2m}
%\sum_{l=0}^m\frac{\abs{\mathcal{P}_\ell^{(\alpha,\beta)}(x)}^2}{\norm{\mathcal{P}_\ell^{(\alpha, \beta)}(x)}^2}
%= \frac{2^{(\alpha+\beta)}}{\pi(1-x)^\frac{2\alpha+1}{2}(1+x)^\frac{2\beta+1}{2}}.
%\end{equation*}
%\end{thm}
By using Theorem \ref{T:main} and the results from the previous section
we will obtain an identity for the Jacobi polynomials.
%\begin{thm}\label{T:jacobimain}
%For nonnegative integers $p, q$,
%\begin{equation}\label{E:identity}
%\lim_{m\to\infty} \frac{1}{m}
%\sum_{l=0}^m\frac{(2\ell+\frac{p+2q}{2})\Gamma(\ell+1)\Gamma(\ell+\frac{p+2q}{2})}
%{\Gamma(\ell+\frac{p+q+1}{2})\Gamma(\ell+\frac{q+1}{2})}
%\abs{\mathcal{P}_\ell^{\(\frac{p+q-1}{2},\frac{q-1}{2}\)}(x)}^2
%= \frac{2^{\frac{p+2q}{2}}}{\pi(1-x)^{\frac{p+q}{2}}(1+x)^\frac{q}{2}}..
%\end{equation}
%\end{thm}
\begin{thm}\label{T:jacobimain}
Let $\alpha=\frac{p+q-1}{2}$ and $\beta=\frac{q-1}{2}$.
Then for $p$ and $q$ given in Table \ref{table:parameters}, we have
\begin{equation}\label{E:identity}
\begin{gathered}
\lim_{m\to\infty} \frac{1}{m}
\sum_{l=0}^m\frac{(2\ell+\alpha+\beta+1)\Gamma(\ell+1)\Gamma(\ell+\alpha+\beta+1)}
{\Gamma(\ell+\alpha+1)\Gamma(\ell+\beta+1)}
\abs{\mathcal{P}_\ell^{\(\alpha,\beta\)}(x)}^2\\
= \frac{2^{\alpha+\beta+1}}{\pi(1-x)^{\frac{2\alpha+1}{2}}(1+x)^\frac{2\beta+1}{2}}.
\end{gathered}
\end{equation}
\end{thm}
\begin{proof}
Let $N=\{u\in M \st d(e,u)=r\}$ for some $r<L$.
Then $N$ is a smooth submanifold of codimension $1$, since the injectivity radius of $M$ is $L$
(see \cite[Chapter \MakeUppercase{\romannumeral 9},Theorem 5.4]{helgason(dgss)}).
Let $\tau=\nu$ be the Riemannian measure on $N$ (i.e. $\psi=1$). Then
$\nu(N)=\frac{2\pi^\frac{p+q+1}{2}}{\Gamma(\frac{p+q+1}{2})\omega^{p+q}}\sin^{p+q}\omega r\cos^q\omega r$.
Therefore according to Theorem \ref{T:main},
\begin{equation}\label{E:sph-harm}
\begin{aligned}
\sum_{\gamma_j < T} \abs{\hat{\tau}(j)}^2
=&\; \sum_{\gamma_j < T} \abs{\<\tau, \Phi_j\>}^2
\sim&\; \frac{\nu(N)}{\pi} T^{1/2} .
\end{aligned}
\end{equation}

%Since $\tau$ is $K$-invariant, $\<\tau, Y_\ell^i\>=0$ if $i\ne 0$ (see Lemma \ref{L:ip}).
%Therefore we have
%\begin{equation*}\label{E:sph-harm-2}
%\begin{aligned}
%%\sum_{\lambda_\ell<T} \abs{\<\tau, Y_\ell^0\>}^2  \sim \frac{\nu(N)}{\pi} T^{1/2}.
%\sum_{\lambda_\ell<T} \abs{\int_M Y_\ell^0(u)~d\nu}^2 
%&\;\sim \frac{\nu(N)}{\pi} T^{1/2}\\
%\sum_{\lambda_\ell<T} \abs{\int_M \phi_\ell(u)~d\nu}^2 
%&\;\sim \frac{\nu(N)}{\pi} T^{1/2}\\
%\sum_{\lambda_\ell<T}
%\abs{\int_M c_\ell\mathcal{P}_\ell^{\(\frac{p+q-1}{2},\frac{q-1}{2}\)}(\cos2\omega\rho(u,e))~d\nu}^2
%&\;\sim \frac{\nu(N)}{\pi} T^{1/2}\\
%\sum_{\lambda_\ell<T} c_\ell^2\abs{\mathcal{P}_\ell^{\(\frac{p+q-1}{2},\frac{q-1}{2}\)}(x)}^2\nu(N)^2
%&\;\sim \frac{\nu(N)}{\pi} T^{1/2}\\
%\sum_{\lambda_\ell<T}c_\ell^2\abs{\mathcal{P}_\ell^{\(\frac{p+q-1}{2},\frac{q-1}{2}\)}(x)}^2
%&\;\sim\frac{T^{1/2}}{\pi\nu(N)}\\
%\sum_{\lambda_\ell<T}c_\ell^2\abs{\mathcal{P}_\ell^{\(\frac{p+q-1}{2},\frac{q-1}{2}\)}(x)}^2
%&\;\sim\frac{T^{1/2}}{\pi}\frac{\Gamma(\frac{p+q+1}{2})
%\omega^{p+q}2^{\frac{p+2q-2}{2}}}{\pi^\frac{p+q+1}{2}(1-x)^\frac{p+q}{2}(1+x)^\frac{q}{2}}.
%\end{aligned}
%\end{equation*}
Therefore
\begin{equation*}
\begin{aligned}
\sum_{\lambda_\ell<T}c_\ell^2\abs{\mathcal{P}_\ell^{\(\frac{p+q-1}{2},\frac{q-1}{2}\)}(x)}^2
&\;= \frac{1}{\nu(N)^2}\sum_{\lambda_\ell<T}
\abs{\int_M c_\ell\mathcal{P}_\ell^{\(\frac{p+q-1}{2},\frac{q-1}{2}\)}(\cos2\omega\rho(u,e))~d\nu}^2\\
&\;= \frac{1}{\nu(N)^2}\sum_{\lambda_\ell<T} \abs{\int_M \phi_\ell(u)~d\nu}^2\\
&\;= \frac{1}{\nu(N)^2}\sum_{\gamma_j<T} \abs{\int_M \Phi_j(u)~d\nu}^2 \quad \text{(By Lemma \ref{L:ip})}\\
&\; \sim\frac{T^{1/2}}{\pi\nu(N)} \quad \text{(By Equation(\ref{E:sph-harm}))}.
\end{aligned}
\end{equation*}

It follows from Lemma \ref{L:projection} that
\begin{equation*}
\lim_{m\to\infty} \frac{1}{m}
\sum_{l=0}^m\frac{(\frac{4\ell+p+2q}{2})\Gamma(\ell+1)\Gamma(\frac{2\ell+p+2q}{2})}
{\Gamma(\frac{2\ell+p+q+1}{2})\Gamma(\frac{2\ell+q+1}{2})}
\abs{\mathcal{P}_\ell^{\(\frac{p+q-1}{2},\frac{q-1}{2}\)}(x)}^2
= \frac{2^{\frac{p+2q}{2}}}{\pi(1-x)^{\frac{p+q}{2}}(1+x)^\frac{q}{2}}.
\end{equation*}
By taking $p=2(\alpha-\beta)$ and $q=2\beta+1$, we get the identity stated in the theorem.
\end{proof}

\section{Some identities for the particular case $(x=-1)$}
Suppose $M$ is a compact rank-one symmetric space other than a sphere.
Let $N=\{u\in M \st d(e,u)=\frac{\pi}{2}\}$.
Then $N$ is the cut locus of $e$. It is well known that $N$ is a smooth manifold (see \cite[Proposition 5.1]{helgasonradon}).
Let $k$ be the codimension of $N$. For values of $k$ and $\nu(N)$ see Table {\ref{table:codimension}}. Then
\begin{cor}
For $p$ and $q$ given in Table \ref{table:parameters}  and $k$ given in Table \ref{table:codimension}, we have
\begin{equation*}
\lim_{m\to\infty} \frac{1}{m^k}
\sum_{l=0}^m\frac{(\frac{4\ell+p+2q}{2})\Gamma(\frac{2\ell+p+2q}{2})\Gamma(\frac{2\ell+q+1}{2})}
{\Gamma(\frac{2\ell+p+q+1}{2})\Gamma(\ell+1)}
=\frac{2}{k}.
\end{equation*}
\end{cor}
\begin{proof}
We have
\begin{equation*}
\begin{aligned}
\sum_{\lambda_\ell<T}c_\ell^2\abs{\mathcal{P}_\ell^{\(\frac{p+q-1}{2},\frac{q-1}{2}\)}(x)}^2
&\;= \frac{1}{\nu(N)^2}\sum_{\lambda_\ell<T}
\abs{\int_M c_\ell\mathcal{P}_\ell^{\(\frac{p+q-1}{2},\frac{q-1}{2}\)}(\cos2\omega\rho(u,e))~d\nu}^2\\
&\;= \frac{1}{\nu(N)^2}\sum_{\lambda_\ell<T} \abs{\int_M \phi_\ell(u)~d\nu}^2\\
&\;= \frac{1}{\nu(N)^2}\sum_{\gamma_j<T} \abs{\int_M \Phi_j(u)~d\nu}^2 \quad \text{(By Lemma \ref{L:ip})}\\
&\; \sim\frac{T^{k/2}}{(4\pi)^{k/2}\Gamma(\frac{k}{2}+1)\nu(N)} \quad \text{(By Theorem \ref{T:main})}.
\end{aligned}
\end{equation*}
Note that
$\mathcal{P}_\ell^{(\alpha,\beta)}(-1)=(-1)^\ell \begin{pmatrix} \ell+\beta \\ \ell \end{pmatrix}$
(see \cite[Equation (4.1.4)]{szego}).

Therefore for $x=-1$, we have
\begin{equation*}
\sum_{\lambda_\ell<T}c_\ell^2\abs{(-1)^\ell \begin{pmatrix} \frac{2\ell+q-1}{2} \\ \ell \end{pmatrix}}^2
\sim \frac{T^{k/2}}{(4\pi)^{k/2}\Gamma(\frac{k}{2}+1)\nu(N)}
\end{equation*}
or equivalently (by using Lemma \ref{L:projection})
%\begin{equation}\label{E:particular}
%\lim_{m\to\infty} \frac{1}{\lambda_m^{k/2}} \sum_{l=0}^m
%\frac{\abs{\begin{pmatrix} \ell+(\frac{q-1}{2}) \\ \ell \end{pmatrix}}^2}
%{\norm{\mathcal{P}_\ell^{\(\frac{p+q-1}{2},\frac{q-1}{2}\)}}^2}
%=\frac{1}{c\(2^\frac{2k+p+2q}{2}\Gamma(\frac{k}{2}+1)\)},
%\end{equation}
%where $c$ is a constant (see Table \ref{table:constant}).

\begin{equation*}
\lim_{m\to\infty} \frac{1}{m^k}
\sum_{l=0}^m\frac{(\frac{4\ell+p+2q}{2})\Gamma(\frac{2\ell+p+2q}{2})\Gamma(\frac{2\ell+q+1}{2})}
{\Gamma(\frac{2\ell+p+q+1}{2})\Gamma(\ell+1)}
=\frac{2}{k}.
\end{equation*}
\end{proof}
\medskip

\begin{center}
\begin{table}[h!]
\caption{Values of $k$ and $\nu(N)$}\label{table:codimension}
\begin{tabular}{|c|c|c|}
\hline
%\toprule
$M$                                               & $k$     & $\nu(N)$\\
%\toprule
\hline
$\C\mathbb{P}^n$                        & $2$      & $\frac{\pi^{n-1}}{\Gamma{(n)}}$\\
\hline
$\mathbb{H}\mathbb{P}^n$          & $4$      & $\frac{\pi^{2(n-1)}}{\Gamma{(2n)}}$\\
\hline
$Ca\mathbb{P}^2$                        & $8$      & $\frac{\pi^4\Gamma(4)}{\Gamma(8)}$\\
\hline
%\bottomrule
\end{tabular}
\end{table}
\end{center}

\section{Acknowledgements}
I sincerely thank my supervisor, Prof. M.~K.~Vemuri, for their guidance,
support, and valuable feedback throughout this research.

%\section{Declarations}
%\subsection{Funding}
%The author was supported by a fellowship from CSIR (File Number :\\09/1217(0077)/2019-EMR-I).
%\subsection{Competing Interests}
%The author has no competing interests to declare that are relevant to the content of this article. 

\bibliographystyle{amsplain}
\bibliography{apij}
\end{document}